\documentclass[1p]{elsarticle}
\usepackage{lineno,hyperref,bm,amsmath,amssymb,algorithm}
\usepackage[]{algpseudocode}
\modulolinenumbers[5]
\newdefinition{rmk}{Remark}

\newcommand\norm[1]{\left\lVert#1\right\rVert}
\newcommand\abs[1]{\left\lvert#1\right\rvert}
\usepackage{subcaption}
\usepackage{caption}
\usepackage{amssymb}

\usepackage{lipsum}
\makeatletter
\def\ps@pprintTitle{%
 \let\@oddhead\@empty
 \let\@evenhead\@empty
 \def\@oddfoot{}%
 \let\@evenfoot\@oddfoot}
\makeatother

\begin{document}

\begin{frontmatter}

\title{Reduced order modelling of nonaffine problems on parameterized NURBS multipatch geometries}


\author[mymainaddress]{Margarita Chasapi\corref{mycorrespondingauthor}}
\cortext[mycorrespondingauthor]{Corresponding author}
\ead{margarita.chasapi@epfl.ch}

\author[mymainaddress]{Pablo Antolin}
\author[mymainaddress,mysecondaryaddress]{Annalisa Buffa}

\address[mymainaddress]{Institute of Mathematics, \'Ecole Polytechnique F\'ed\'erale de Lausanne, Lausanne, Switzerland }
\address[mysecondaryaddress]{Instituto di Matematica Applicata e Tecnologie Informatiche 'E. Magenes' (CNR), Pavia, Italy}

\begin{abstract}
This contribution explores the combined capabilities of reduced basis methods and IsoGeometric Analysis (IGA) in the context of parameterized partial differential equations. The introduction of IGA enables a unified simulation framework based on a single geometry representation for both design and analysis. The coupling of reduced basis methods with IGA has been motivated in particular by their combined capabilities for geometric design and solution of parameterized geometries. In most IGA applications, the geometry is modelled by multiple patches with different physical or geometrical parameters. In particular, we are interested in nonaffine problems characterized by a high-dimensional parameter space. We consider the Empirical Interpolation Method (EIM) to recover an affine parametric dependence and combine domain decomposition to reduce the dimensionality. We couple spline patches in a parameterized setting, where multiple evaluations are performed for a given set of geometrical parameters, and employ the Static Condensation Reduced Basis Element (SCRBE) method. At the common interface between adjacent patches a static condensation procedure is employed, whereas in the interior a reduced basis approximation enables an efficient offline/online decomposition. The full order model over which we setup the RB formulation is based on NURBS approximation, whereas the reduced basis construction relies on techniques such as the Greedy algorithm or proper orthogonal decomposition (POD). We demonstrate the developed procedure using an illustrative model problem on a three-dimensional geometry featuring a multi-dimensional geometrical parameterization.
\end{abstract}


\end{frontmatter}

\section{Introduction}
Isogeometric analysis was introduced as a paradigm that unifies the design and analysis process by employing the same representation for the geometry and approximation of the solution \cite{Hughes2005}. In this work we combine isogeometric analysis (IGA) and reduced basis (RB) methods for the solution of partial differential equations (PDEs) on parameterized geometries. Our motivation is to exploit the capabilities of splines for geometric representation in the context of reduced basis solution of parameterized PDEs. In the past, there have been several successful applications of reduced basis methods in isogeometric analysis. Nonaffine geometrical parameterizations were introduced in \cite{Manzoni2015,Salmoiraghi2016} in the context of fluid dynamics problems. A certified Greedy-based reduced basis method was proposed for affinely parametric problems approximated by NURBS in \cite{Devaud2017}. Isogeometric analysis was further combined with proper orthogonal decomposition (POD) for parabolic problems in \cite{Zhu2017} as well as with interpolation \cite{Garotta2020} and deep learning \cite{Fresca2017} to allow non-intrusive reduced order modelling. Most of these works focus on problems with a few parameters.

In many practical applications, geometries are often represented by multiple patches. These can be viewed as subdomains with different physical or geometrical parameters. The dimension of the latter might be quite high with increasing geometric complexity. Isogeometric analysis was also used in the past to construct reduced order models for complex geometries based on the POD technique and the idea of isotopological meshing in \cite{Maquart2020}. The new aspect of our work is to exploit domain decomposition within an efficient RB framework for nonaffine geometrical parameterizations that may be characterized by a large number of parameters. The application of the RB method to such problems is not straightforward. To tackle nonaffine problems, we consider the Empirical Interpolation Method (EIM) \cite{Barrault2004}. This technique recovers affine dependence and enables an efficient offline/online decomposition. Moreover, the efficiency of RB methods depends highly on the number of parameters. The Static Condensation Reduced Basis Element Method (SCRBE) was introduced to overcome such shortcomings and allow to solve problems of industrial relevance within a reduced basis framework \cite{Huynh2013a,Huynh2013b,Vallaghe2015}. In the context of SCRBE, the RB approximation is combined with a static condensation procedure to facilitate efficient reduction in the interior of components and at interfaces. The method was further extended to enable additional reduction at interfaces \cite{Eftang2013,Eftang2013b,Smetana2016} and nonlinearities \cite{Ballani2018}.

This contribution is structured as follows: Section 2 presents the parameterized model problem that we use for illustration purposes throughout this work. In Section 3 we briefly review the main concepts related to splines and multipatch geometries. An overview of the reduced basis method and space construction is provided in Sections 4 and 5, respectively. In Section 6 we move on to the Empirical Interpolation Method and recast its formulation in the context of isogeometric analysis. Section 7 gives a brief overview of the SCRBE method. The presented procedure is applied to a numerical example in Section 8. Finally, we summarize the main conclusions that can be drawn from this study and provide some outlook for future work.

\section{Parameterized model problem}\label{sec:model}
Throughout this work, we consider the setting of linear elliptic parameterized partial differential equations. In the following, the parameterized Poisson equation will serve as model problem for ease of exposition.

Let us consider a parameterized domain $\Omega(\bm{\mu}) \subset \mathbb{R}^d$, where $d$ is the dimension of the physical space of the problem at hand. The domain is described by the geometrical parameters $\bm{\mu} \in \mathcal{P} \subset \mathbb{R}^P$, where $\mathcal{P}$ is the parameter space and $P$ is the number of parameters. The continuous formulation of the problem reads in strong form: for any $\bm{\mu} \in \mathcal{P}$, find $u \in H_{0,{\Gamma_D}}^1(\Omega(\bm{\mu}))$ such that
\begin{equation}\label{eq1}
\begin{cases}
- \Delta u &= \tilde{f} \ \text{in} \ \Omega({\bm{\mu}}) \\
 \quad \ \ \displaystyle  u &= 0 \ \text{on} \ \Gamma_D({\bm{\mu}}) \\
 \quad \ \displaystyle  \frac{\partial{u}}{\partial{\boldsymbol{n}}} &= 0 \ \text{on} \ \Gamma_{N}({\bm{\mu}}), 
\end{cases}
\end{equation}
where $\Gamma_D({\bm{\mu}})$ and $\Gamma_N({\bm{\mu}})$ denote the Dirichlet and Neumann part of the boundary, accordingly, while it holds $\overline{\Gamma}_D({\bm{\mu}})\cup\overline{\Gamma}_N({\bm{\mu}}) = \partial{\Omega}({\bm{\mu}})$ and $\Gamma_D({\bm{\mu}}) \cap \Gamma_N({\bm{\mu}}) = \varnothing$. We remark that for the scalar case at hand, $ H_{0,{\Gamma_D}}^1(\Omega(\bm{\mu})) \subset H^1(\Omega(\bm{\mu}))$ is the subspace of functions in $H^1(\Omega(\bm{\mu}))$ with vanishing traces on the boundary $\Gamma_D$. Moreover, $\tilde{f} \in L^2(\Omega(\bm{\mu}))$ represents the source term and  $\boldsymbol{n}$ is the outward normal to $\partial{\Omega}({\bm{\mu}})$. For the sake of simplicity, we assume homogeneous Dirichlet and Neumann boundary conditions without loss of generality.
We can express the discrete weak formulation of the parameterized problem in Equation \ref{eq1} as: find $u_h \in V_h$ such that 
\begin{equation}\label{eq2}
\alpha(u_h,v_h;\bm{\mu}) = f(v_h;\bm{\mu}), \qquad \forall v_h \in V_h,
\end{equation}
where $V_h \subset H_{0,\Gamma_D}^1(\Omega(\bm{\mu}))$ is a finite-dimensional subspace. 
The associated parameterized bilinear form $\alpha(\cdot,\cdot;\bm{\mu})$ and the linear functional $f(\cdot;\bm{\mu})$ read:
\begin{equation}\label{eq3}
\begin{aligned}
\alpha(u_h,v_h;\bm{\mu}) &= \int_{{\Omega}(\bm{\mu})} \nabla{u}_h \cdot \nabla{v}_h \, \textrm{d}\Omega, \\
f(v_h;\bm{\mu}) &= \int_{{\Omega}(\bm{\mu})} \tilde{f} v_h \, \textrm{d}\Omega.
\end{aligned}
\end{equation}
The discrete approximation leads to the following parameterized linear system of dimension $\mathcal{N}_h = \text{dim}(V_h)$
\begin{equation}\label{eq4}
{\bf{A}}(\bm{\mu}){\bf{u}}_h(\bm{\mu}) = {\bf{f}}(\bm{\mu}),
\end{equation}
where ${\bf{A}}({\bm{\mu}})$ $\in \mathbb{R}^{\mathcal{N}_h\times \mathcal{N}_h}$ is the stiffness matrix corresponding to the differential operator, ${\bf{f}}({\bm{\mu}}) \in \mathbb{R}^{\mathcal{N}_h}$ is the vector representing the source term and $\mathcal{N}_h$ is the number of degrees of freedom. 

\section{B-splines and parameterized multipatch geometries}\label{sec:iga}
In this section, we briefly review the concept of B-splines and multipatch geometries. A more detailed exposition is available in \cite{Hughes2005,Cottrell2009,Piegl1995}. 
Let us first introduce the parametric domain $\hat{\Omega} = [0,1]^{\hat{d}}$, which is independent of the geometrical parameters and represents the counterpart of our parameter-dependent domain $\Omega(\bm{\mu})$ defined previously in the physical space. We assume that the dimension is the same for both the physical and parametric domain $d=\hat{d}$ without loss of generality. 
Then, we introduce a univariate B-spline basis function $b_{i_j,p_j}^j$ and a \emph{knot vector} $\Xi = \{\xi_1,...,\xi_{n+p+1}\}$ in $\hat{\Omega}$. The integers $p$ and $n$ denote the degree  and the number of basis functions, accordingly. The definition of B-splines $\mathcal{B}_{\boldsymbol{i},\boldsymbol{p}}(\boldsymbol{\xi})$ based on the tensor product of univariate B-spline basis functions reads:
\begin{equation}\label{eq5}
\mathcal{B}_{\bf{i},\bf{p}}(\boldsymbol{\xi}) = \prod_{j=1}^{{d}} b_{i_j,p_j}^j(\xi^j),
\end{equation}
where ${d} \in \{1,2,3\}$ is the corresponding dimension. Moreover, the multi-index ${\bf{i}}=(i_1,..,i_{{d}})$ denotes the position in the tensor-product structure and ${\bf{p}}=(p_1,..,p_{{d}})$ are the polynomial degrees corresponding to the parametric directions $\boldsymbol{\xi}=(\xi^1,..,\xi^{{d}})$. In the following, we will assume that the vector $\bf{p}$ is identical in all parametric directions and can be replaced by a scalar value $p$. Given a multiplicity $k$ at internal knots, the continuity of the B-spline basis is $C^{p-k}$ at every internal knot and $C^{\infty}$ elsewhere. This concept can be extended to rational B-splines (NURBS) in a straightforward manner \cite{Piegl1995}. 

\begin{figure}[!h]
     \centering
     \begin{subfigure}[b]{0.32\textwidth}
         \centering
         \includegraphics[width=\textwidth]{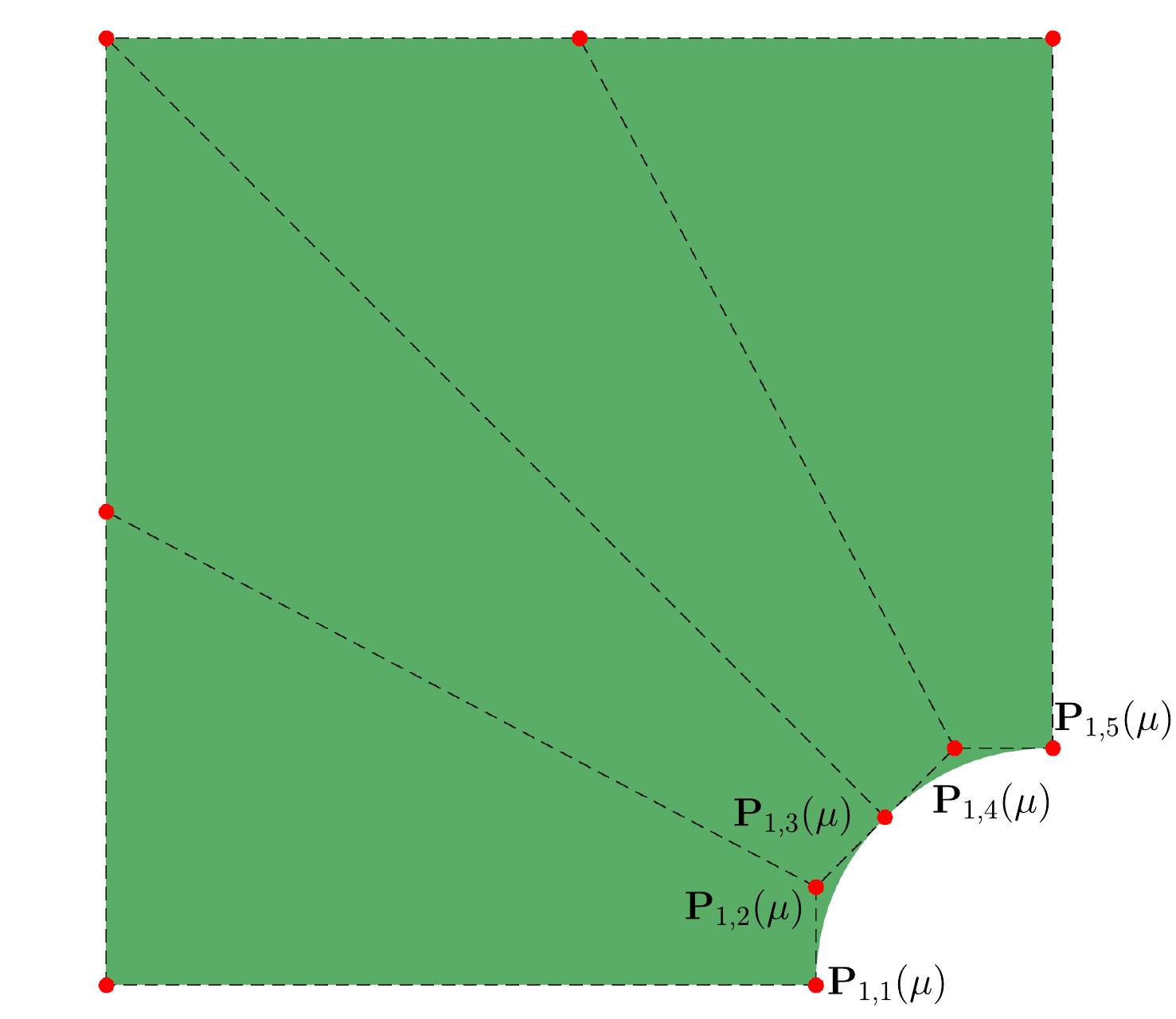}
         \caption{$\mu = 1$}
         \label{fig:CP1}
     \end{subfigure}
     \begin{subfigure}[b]{0.32\textwidth}
         \centering
         \includegraphics[width=\textwidth]{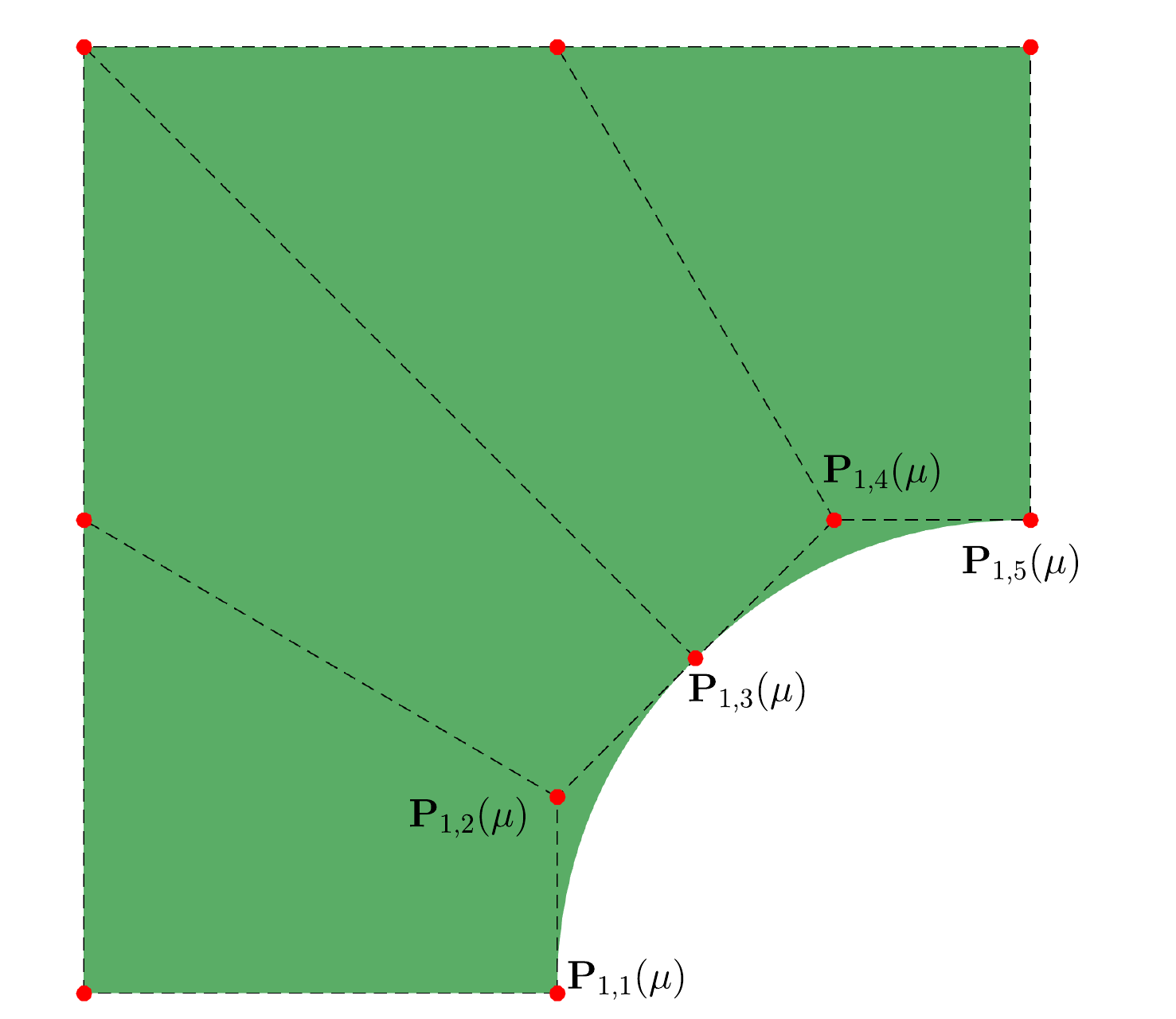}
         \caption{$\mu =2$}
         \label{fig:CP2}
     \end{subfigure}
     \begin{subfigure}[b]{0.32\textwidth}
         \centering
         \includegraphics[width=\textwidth]{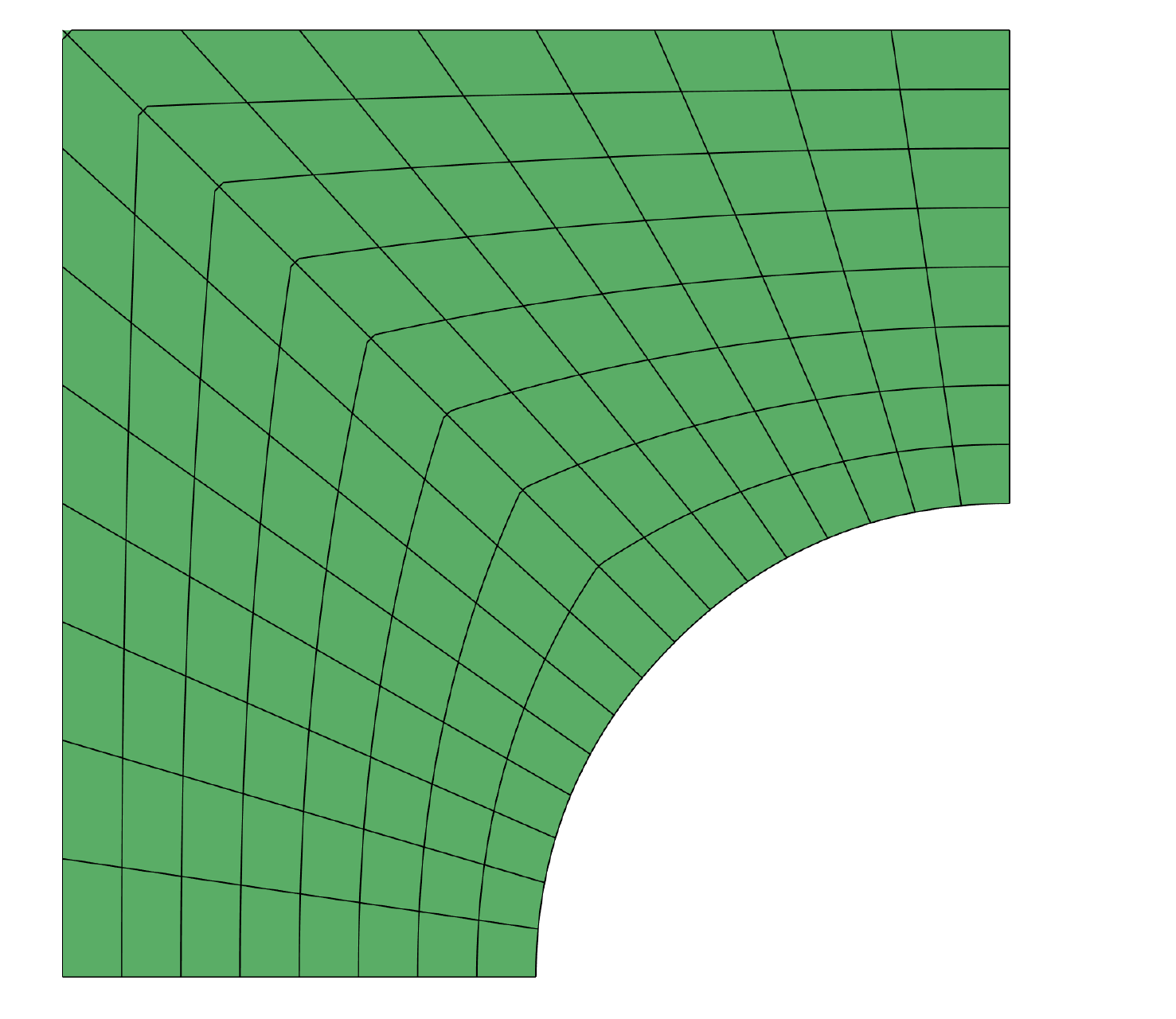}
         \caption{$\mu = 2$}
         \label{fig:CP_mesh}
     \end{subfigure}
        \caption{Coarse geometries (a) and (b) for different values of a geometrical parameter $\mu$ representing the radius of the circular hole. The refined mesh (c) is employed for the analysis.}
      \label{fig:CP_mu}
\end{figure}

Having the definition of B-splines at hand, we can now introduce the geometric map ${\bf{F}} : \hat{\Omega} \to \Omega$ as
\begin{equation}\label{eq6}
{\bf{F}}(\boldsymbol{\xi};\bm{\mu}) = \sum_{\bf{i}} \mathcal{B}_{{\bf{i}},p}(\boldsymbol{\xi})\bf{P_i}(\bm{\mu}),
\end{equation}
where the parameterized control points are denoted by $\bf{P_i}(\bm{\mu})$. We assume that ${\bf{F}}$ is smooth and invertible with piecewise smooth inverse. Moreover, in the case of NURBS we assume that the weights are parameter-independent. It is worth-wile remarking that we can deform the physical geometry by moving the control points $\bf{P_i}(\bm{\mu})$ for a given value of the geometrical parameters $\bm{\mu}$ as illustrated in Figure \ref{fig:CP_mu}. If the  number of parameterized control points is large, the number of geometrical parameters $P$ might grow drastically. Therefore it seems advantageous to define the latter on a very coarse geometry and then perform refinement in the analysis phase. This is possible in isogeometric analysis since the geometry is preserved during the refinement process. The reader is referred to \cite{Cottrell2007} for more details on refinement strategies. A similar concept is commonly adopted in isogeometric shape optimisation, which makes this choice of shape parameterization easily adaptable to parametric optimization studies \cite{Nagy2013,Kiendl2014}. We note that even if the geometry changes, we can obtain a low-dimensional parameterization by introducing a mapping between the exact geometry and the one employed to define the parameters.

\begin{figure}[!htb]
	\centering
	\includegraphics[width=0.7\textwidth]{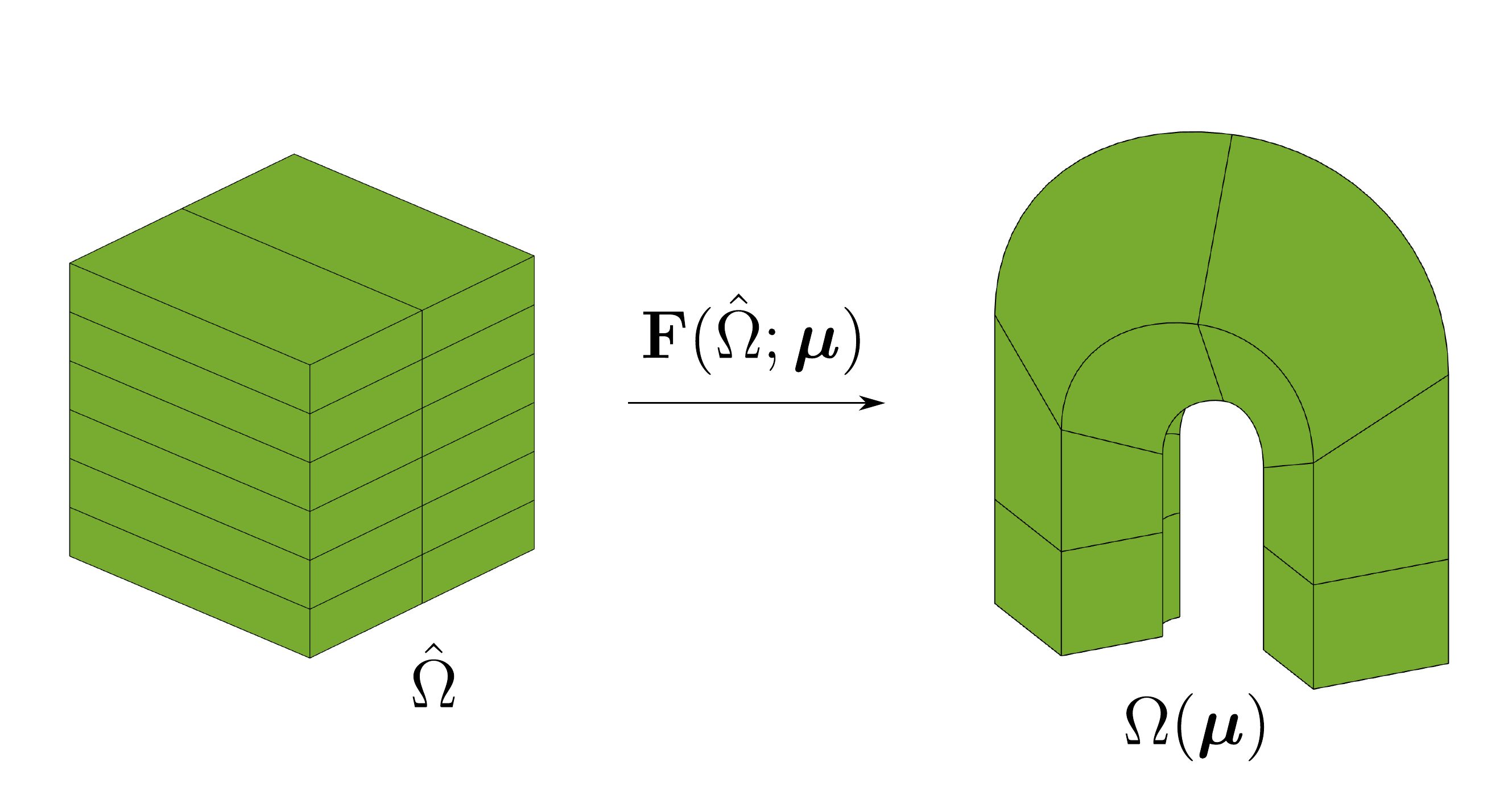} \\
	\caption{Geometry of horseshoe obtained by B-spline mapping}\label{fig:mapping}
\end{figure}

With the definition of the geometric map at hand, the physical domain can be expressed as
\begin{equation}\label{eq7}
\Omega(\bm{\mu}) = {\bf{F}}(\hat{\Omega};\bm{\mu}). 
\end{equation}
The spline mapping is illustrated for an exemplary geometry in Figure \ref{fig:mapping}. The B-spline space on the physical domain is defined as:
\begin{equation}\label{eq8}
S_h({{{{\Omega}}({\bm{\mu}})}}) = \text{span}\{\mathcal{B}_{{\bf{i}},p} \circ {\bf{F}}^{-1} \}. 
\end{equation}

So far we have presented the case of a single patch, that is the computational domain is obtained with a single geometric map from the parametric domain to the physical one. In many practical applications, the geometry of interest is complex and needs to be represented by multiple patches. These can be understood as subdomains that can be characterized by different physical properties. In the following we assume $C^0$ continuity between patches, although in principle this can be increased \cite{Cottrell2007,Kapl2017,Kapl2017b}. 

Let us consider that the parameterized domain $\Omega(\bm{\mu})$ is obtained as a union of $n_p$ patches such that
\begin{equation}\label{eq9}
\overline{\Omega}(\bm{\mu}) = \bigcup_{k=1}^{n_p} {\overline{\Omega}}_k(\bm{\mu}^{(k)}),
\end{equation}
where $\Omega_k(\bm{\mu}^{(k)}) \cap \Omega_l(\bm{\mu}^{(l)}) = \varnothing$ for $l \ne k$. Here $\bm{\mu}^{(k)} \in \mathcal{P}_k \subset \mathbb{R}^{P_k}$ are the parameters associated to the $k$-th patch with $\mathcal{P} \subseteq \oplus_{k=1}^{n
_p}\mathcal{P}_k$ and $\mathcal{P}_k,P_k$ the corresponding parameter space and number of parameters. For each subdomain 
$\Omega_k(\bm{\mu}^{(k)})$, the parametric domain $\hat{\Omega} = [0,1]^{\hat{d}}$ is the same and it holds
\begin{equation}\label{eq10}
\Omega_k(\bm{\mu}^{(k)}) = {\bf{F}}_k(\hat{\Omega};\bm{\mu}^{(k)}),
\end{equation}
where ${\bf{F}}_k$ denotes the geometric map associated to the $k$-th patch. We can now define the approximation space $S_h^{(k)}$ in the patch $\Omega_k(\bm{\mu}^{(k)})$ as 
\begin{equation}\label{eq11}
S_h^{(k)}({{{{\Omega}}_k({\bm{\mu}}^{(k)})}}) = \text{span}\{\mathcal{B}_{{\bf{i}},p}^{(k)} \circ {\bf{F}}_k^{-1} \},
\end{equation}
where $\mathcal{B}_{{\bf{i}},p}^{(k)}$ are the B-splines corresponding to the $k$-th patch. In the following we will assume that patches are matching at the common interfaces, that is for $\Gamma_{k,l}=\partial{\Omega_k}\cap\partial{\Omega_l} \ne \varnothing$ for all $\beta^{(k)} \in {\mathcal{B}_{{\bf{i}},p}^{(k)}}$, there exists a $\beta^{(l)} \in {\mathcal{B}_{{\bf{i}},p}^{(l)}}$ 
such that ${\beta^{(k)}}{\vert{_{\Gamma_{k,l}}}} = {\beta^{(l)}}{\vert_{\Gamma_{k,l}}}$, and vice versa. The space $H_{0,{\Gamma_D}}^1(\Omega(\bm{\mu}))$ is then approximated by the spline space:
\begin{equation}\label{eq11b}
V_h = \{ u  \in H_{0,{\Gamma_D}}^1(\Omega(\bm{\mu})) : u \vert_{\Omega_k} \in S_h^{(k)}({{{{\Omega}}_k({\bm{\mu}}^{(k)})}}) \}.
\end{equation}

\section{Reduced basis method for parameterized PDEs}\label{sec:rb}
In the following, we are interested in solving Equation \eqref{eq2} for several values of the parameters $\bm{\mu}$. In a real-time and many-query context, a full order approximation of the solution for each new parameter would be computationally inefficient. To reduce the computational cost of approximating the solution for every input parameter, we consider the reduced basis method.

Let us first introduce the discrete solution manifold
\begin{equation}\label{eq12}
\mathcal{M}_h = \{u_h(\bm{\mu}) \in V_h : \bm{\mu} \in \mathcal{P} \} \subset V_h,
\end{equation}
where $u_h(\bm{\mu})$ is the solution of the full order problem in Equation \eqref{eq2}. In the following we consider $\mathcal{N}_h$ to be sufficiently large so that the approximation of the solution $u_h(\bm{\mu})$ can be seen as a \emph{high-fidelity} or \emph{truth} solution. The RB method seeks a low-dimensional approximation of $\mathcal{M}_h$ in order to approximate the solution $u_h(\bm{\mu})$ as a linear combination of suitable basis functions based on a Galerkin projection onto a subspace $V_N \subset V_h$ of dimension $N \ll \mathcal{N}_h$. The dimension $N$ should be as small as possible to reduce the computational cost. A reduced basis of $V_N$ can be constructed based, for example, on a Greedy algorithm or Proper Orthogonal Decomposition. These techniques will be further addressed in Section \ref{sec:basis}.
Considering a Galerkin projection, the reduced basis approximation of the problem in Equation \eqref{eq2} reads: find $u_N \in V_N$ such that
\begin{equation}\label{eq13}
\alpha(u_N,v_N;\bm{\mu}) = f(v_N;\bm{\mu}), \qquad \forall v_N \in V_N.
\end{equation}
The linear system of Equation \ref{eq4} is rewritten as
\begin{equation}\label{eq13b}
{\bf{A}}_N(\bm{\mu}) {\bf{u}}_N(\bm{\mu}) = {\bf{f}}_N(\bm{\mu}).
\end{equation}
One crucial assumption concerning the efficiency of the RB method is that the parameterized bilinear form $\alpha(\cdot,\cdot;\bm{\mu})$ and the linear functional $f(\cdot;\bm{\mu})$ admit an affine decomposition in the form of
\begin{equation}\label{eq14}
\alpha(u_h,v_h;\bm{\mu}) =  \sum_{q=1}^{Q_{\alpha}} \theta_q^{\alpha} (\bm{\mu}) \alpha_q(u_h,v_h), \quad f(v_h;\bm{\mu})=\sum_{q=1}^{Q_f} \theta_q^f(\bm{\mu}) f_q(u_h,v_h),
\end{equation}
where $\theta_q^{\alpha} : \mathcal{P} \to \mathbb{R} $ for $q=1,..,Q_{\alpha}$ and $\theta_q^f : \mathcal{P} \to \mathbb{R} $ for $q=1,..,Q_f$ are $\bm{\mu}$-dependent parameter functions, whereas $\alpha_q : V_h \times V_h \to \mathbb{R}$ and  $f_q : V_h \to \mathbb{R}$ are ${\bm{\mu}}$-independent forms. The latter are independent of the input parameter ${\bm{\mu}}$, therefore they are assembled and stored in the offline phase. In the online phase, the assembly is inexpensive considering the affine parametric dependence in Equation \eqref{eq14} for $Q_{\alpha}, Q_f \ll \mathcal{N}_h$. In particular, during the offline phase the parameter-independent forms $\alpha_q, f_q$ are pre-computed and their reduced counterparts are obtained once a reduced basis for $V_N$ is constructed. In the online phase, the assembly requires only the evaluation of the parameter-dependent functions $\theta_q^{\alpha}, \theta_q^f$. We remark that this decomposition is not always fulfilled. For instance, this is the case of the problem with geometrical parameters at hand. In the following, we will consider a hyper-reduction strategy based on the Empirical Interpolation Method \cite{Barrault2004} to recover the affine dependence in respect to the geometrical parameters. A discussion on this procedure in the context of spline approximations is given in Section \ref{sec:eim}.

\section{Reduced basis space construction}\label{sec:basis}
Let us now elaborate on how to obtain a basis for the reduced problem. We will consider for our numerical studies techniques such as the Greedy algorithm and Proper Orthogonal Decomposition (POD). The reader is referred to \cite{QMN_RBspringer,Hesthaven2016} for a more detailed overview of these techniques.

\subsection{Greedy algorithm}\label{sec:greedy}
We first recall the Greedy algorithm to construct iteratively a basis of the subspace $V_N$ according to a suitable optimality criterion. Considering the $n$-th iteration with $1 \le n \le N_{max}$ and a given sample set of parameter values $S_n = \{\bm{\mu}^1,...,\bm{\mu}^n\} \subset \mathcal{P}_{train} \subset \mathcal{P}$, the corresponding subspace reads
\begin{equation}\label{eq15}
V_N=\text{span}\{u_h(\bm{\mu}^1),...,u_h({\bm{\mu}}^n)\}.
\end{equation}
We can now define the snapshot matrix ${\mathbb{S}} \in \mathbb{R}^{N_h \times n}$ as 
\begin{equation}\label{eq16}
{\mathbb{S}}=[{\bf{u}}_h(\bm{\mu}^1),...,{\bf{u}}_h(\bm{\mu}^n)] \in \mathbb{R}^{N_h \times n},
\end{equation}
where ${\bf{u}}_h(\bm{\mu}^i)$ is the solution of the linear system in Equation \ref{eq4} for  $1\le i \le n$. Then we can apply the Gram-Schmidt orthonormalization process to the snapshot matrix ${\mathbb{S}}$ and obtain an orthonormal basis for $V_N$ as 
\begin{equation}\label{eq17}
{\mathbb{V}}=[{\bm{\zeta}}_1,...,{\bm{\zeta}}_n] \in \mathbb{R}^{N_h \times n}.
\end{equation}
\par The sample parameters $\bm{\mu}_n$ are selected in the offline phase from a sufficiently fine training sample set $\mathcal{P}_{train} \in \mathcal{P}$. The main idea behind the Greedy algorithm is to retain at each iteration the snapshot which is worst approximated by the previously computed basis, that is
\begin{equation}\label{eq18}
{\bm{\mu}}^{n+1}=\underset{{\bm{\mu}}\in \mathcal{P}_{train}}{\arg\max} \Delta_n({\bm{\mu}}), 
\end{equation}
where $\Delta_n({\bm{\mu}})$ is an a posteriori error estimator \cite{QMN_RBspringer} such that
\begin{equation}\label{eq18b}
\norm{{\bf{u}}_h(\bm{\mu}) - \mathbb{V} {\bf{u}}_n(\bm{\mu})}_{\mathbb{X}_h} \le {\Delta_n}(\bm{\mu}), \ \forall {\bm{\mu}} \in \mathcal{P}_{train}.
\end{equation}
Here ${\bf{u}}_n(\bm{\mu})$ denotes the solution of the RB problem in Equation \eqref{eq13b} and ${\mathbb{X}_h}$ the matrix associated to the inner product in $H_{0,{\Gamma_D}}^1$. The algorithm is terminated at the $N$-th iteration when the estimator ${\Delta_N}(\bm{\mu})$ is smaller than a prescribed tolerance $\varepsilon$, that is 
\begin{equation}\label{eq19}
\underset{\bm{\mu}\in\mathcal{P}_{train}}{\max} {\Delta_N}(\bm{\mu}) \le \varepsilon,
\end{equation}
whereas $N$ is the final size of the RB space. The evaluation of the error estimator should be independent of the dimension $\mathcal{N}_h$, therefore we compute this considering the affine assumption and projection of the solution onto the low-dimensional space $V_N$. We remark that in this work we consider an inexpensive residual-based a posteriori indicator. The interested reader is referred to the literature for a detailed discussion \cite{Manzoni2015}.

\subsection{Proper orthogonal decomposition}
In the following we consider the construction of the reduced basis based on the Proper Orthogonal Decomposition (POD). This technique relies on the singular value decomposition algorithm (SVD) to reduce the dimension of the original system by extracting a set of orthonormal basis functions \cite{Quarteroni2017}. These can can be viewed as modes retaining most of the energy of the original system.
Let us again consider a fine training sample set $\mathcal{P}_{train} = \{\bm{\mu}^1,...,\bm{\mu}^{N_s}\} \subset \mathcal{P}$ of dimension $N_s$. Then we can construct the snapshots matrix ${\mathbb{S}} \in \mathbb{R}^{\mathcal{N}_{h} \times N_s}$ as
\begin{equation}\label{eq20}
{\mathbb{S}}= [{{\bf{u}}_h(\bm{\mu}^1)},...,{\bf{u}}_h(\bm{\mu}^{N_s})].
\end{equation}
The SVD of the snapshot matrix ${\mathbb{S}}$ reads:
\begin{equation}\label{eq21}
{\mathbb{S}}= \mathbb{U}\boldsymbol{\Sigma}\mathbb{Z}^T,
\end{equation}
where the matrices $\mathbb{U}=[\bm{\zeta}_1,...,\bm{\zeta}_{N_h}] \in \mathbb{R}^{N_h \times N_h}$, $\mathbb{Z}=[\bm{\psi}_1,...,\bm{\psi}_{N_s}] \in \mathbb{R}^{N_s \times N_s}$ are orthogonal with columns containing the left and right singular vectors of ${\mathbb{S}}$, accordingly, and $\boldsymbol{\Sigma}=\text{diag}(\sigma_1,...,\sigma_r)$ with singular values $\sigma_1 \ge \sigma_2 \ge ... \sigma_r$, where $r \le N_s$ being $r$ the rank of $\mathbb{S}$. The POD basis of dimension $N$ is then defined as the set of first $N$ left singular vectors of $\mathbb{S}$ as
\begin{equation}\label{eq22}
\mathbb{V} = [\bm{\zeta}_1,...,\bm{\zeta}_{N}] \in \mathbb{R}^{N_h \times N},
\end{equation}
which correspond to the $N$ largest singular values.
The POD basis is orthonormal by construction. We can select the dimension of the reduced basis $N$ such that the error in the POD basis is smaller than a prescribed tolerance $\varepsilon_{POD}$. The error can be viewed as the sum of the squares of the singular values associated to the neglected POD modes \cite{QMN_RBspringer}. Therefore, we consider $N$ as the smallest integer such that
\begin{equation}\label{eq24}
 1 - \frac{\sum_{i=1}^{N}\sigma_i^2}{\sum_{i=1}^{r}\sigma_i^2} \le \varepsilon_{POD},
\end{equation}
that is the energy captured by the last neglected modes is smaller than or equal to $\varepsilon_{POD}$. 
\par Although the construction of the POD basis is a straightforward task, it requires to compute a priori a sufficiently large amount of high-fidelity solutions in order to obtain a reduced basis with good approximation properties. This drawback can be overcome by constructing the reduced basis with an adaptive technique such as the Greedy algorithm, where only selected snapshots are computed. Throughout this work we will mainly consider the Greedy algorithm to construct RB approximations. In the context of domain decomposition, we will take advantage of the simple construction of the POD basis to construct suitable modes and reduce the dimensionality at common interfaces between patches. A discussion on this aspect will follow in Section \ref{sec:scrbe}.

\section{Empirical interpolation method for spline approximations}\label{sec:eim}
In this work we are interested in the general case of geometrical parameterizations that might entail a nonaffine parametric dependence. At this point we recall that the affine dependence assumption is crucial for the efficiency of the RB method. Therefore, an additional reduction stage is required to recover the affine decomposition of $\bm{\mu}$-dependent forms introduced in Equation \eqref{eq14}. In the following, we will consider the Empirical Interpolation Method (EIM) to ensure an efficient offline/online decomposition. EIM is an interpolation technique, which has been successfully applied to RB methods in order to recover an affine approximation of parameter-dependent functions \cite{Barrault2004}. In this section we will recast its formulation in the context of spline approximations and briefly review the main concepts associated to this hyper-reduction procedure.

\subsection{Weak formulation in reference domain}
Let us first recall the weak formulation of the original problem in Equation \eqref{eq3}. The first step required by the EIM procedure is to pull back the weak formulation to a $\bm{\mu}$-independent reference configuration and apply the EIM algorithm to each term associated to the geometrical parameters. In the context of finite element approximations, it is common practise to define a map from a reference domain to the current configuration. This operation requires to know the analytic expression of the map and its gradient, which is not always an easy task. At this point, we take advantage of the spline definition and the parametric domain $\hat{\Omega} = [0,1]^{\hat{d}}$ introduced in Section \ref{sec:iga}. This allows us to exploit the exact geometric map $\bf{F}$ \eqref{eq6} and define $\hat{\Omega}$ as a $\bm{\mu}$-independent reference domain to apply the EIM procedure \cite{Rinaldi2015}. We can easily obtain the transformed weak formulation by a change of variables as
\begin{equation}\label{eq25}
\begin{aligned}
\alpha(u_h,v_h;\bm{\mu}) &= \int_{{\Omega}(\bm{\mu})} \nabla{u}_h \cdot \nabla{v}_h \, \textrm{d}\Omega =  \int_{\hat{\Omega}} (\nabla\hat{u}_h {\bf{DF}}^{-1}) \cdot (\nabla\hat{v}_h {\bf{DF}}^{-1})|\text{det}{\bf{DF}}| \, \textrm{d}\hat\Omega,\\
f(v_h;\bm{\mu}) &= \int_{{\Omega}(\bm{\mu})} \tilde{f} v_h \, \textrm{d}\Omega = \int_{\hat{\Omega}} \tilde{f}({\bf F}) \ \hat{v}_h |\text{det}{\bf{DF}}| \ \textrm{d}\hat\Omega,
\end{aligned}
\end{equation}
where ${\bf{DF}}$ denotes the Jacobian matrix of the geometric map $\bf{F}$, ${\bf{DF}}^{-1}$ its inverse and $\text{det}{\bf{DF}}$ its determinant. We can now observe that the map ${\bf{F}}(\boldsymbol{\xi};\bm{\mu})$ accounts for the parametric dependency of the integrals. This allows us to define the parameter-dependent family of functions $\mathcal{G}_a,\mathcal{G}_f$ for the EIM procedure as
\begin{equation}\label{eq26}
\begin{aligned}
\mathcal{G}_{\alpha} &= \{g_{\alpha}(\boldsymbol{\xi};\bm{\mu}), \bm{\mu} \in \mathcal{P} \
| \ g_{\alpha}(\boldsymbol{\xi};\bm{\mu}) = 
{\bf{DF}}^{-1}(\boldsymbol{\xi};\bm{\mu}){\bf{DF}}^{-T}(\boldsymbol{\xi};\bm{\mu})|\text{det}{\bf{DF}}(\boldsymbol{\xi};\bm{\mu})|\}, 
\\
\mathcal{G}_f &= \{g_{f}(\boldsymbol{\xi};\bm{\mu}), \bm{\mu} \in \mathcal{P} \ | \ g_{f}(\boldsymbol{\xi};\bm{\mu}) = 
\tilde{f}({\bf{F}}(\boldsymbol{\xi};\bm{\mu}))|\text{det}{\bf{DF}}(\boldsymbol{\xi};\bm{\mu})|\},
\end{aligned}
\end{equation}
for the bilinear form and linear functional, respectively. Here, ${\bf{DF}}^{-T}$ represents the transposed inverse of the Jacobian.
We remark that these functions change if another differential operator is considered, though in principle one can construct them in a similar manner. Now that we have defined the nonaffine functions for the problem at hand we formulate the reduced problem with the EIM approximation.

\subsection{RB problem with EIM approximation}
The main idea behind EIM is to approximate elements of $\mathcal{G}_{\alpha},\mathcal{G}_f$ by constructing an interpolant $\mathcal{I}_M^{\boldsymbol{\xi}}$ for the functions $g_{\alpha}(\boldsymbol{\xi};\bm{\mu})$ and $g_{f}(\boldsymbol{\xi};\bm{\mu})$. We seek an approximate affine expansion of the functions in Equation \eqref{eq26} in the form of 
\begin{equation}\label{eq27}
\begin{aligned}
g_{\alpha}^M(\boldsymbol{\xi};\bm{\mu}) &= \mathcal{I}_M^{\boldsymbol{\xi}} g_{\alpha}(\boldsymbol{\xi};\bm{\mu}) = \sum_{m_1=1}^{M_{\alpha}} \theta_{m_1}^{\alpha}(\bm{\mu})\phi_{m_1}^{\alpha}(\boldsymbol{\xi}),
\\
g_f^M(\boldsymbol{\xi};\bm{\mu}) &= \mathcal{I}_M^{\boldsymbol{\xi}} g_f(\boldsymbol{\xi};\bm{\mu}) = \sum_{m_2=1}^{M_f} \theta_{m_2}^{f}(\bm{\mu})\phi_{m_2}^{f}(\boldsymbol{\xi}),
\end{aligned}
\end{equation}
where $\theta_{m_1}^{\alpha}(\bm{\mu}),\theta_{m_2}^{f}(\bm{\mu})$ are $\bm{\mu}$-dependent functions and $\phi_{m_1}^{\alpha},\phi_{m_2}^{f}$ are $\bm{\mu}$-independent basis functions for $m_1=1,...,M_{\alpha}$ and $m_2=1,...,M_f$, respectively. The superscript $\boldsymbol{\xi}$ denotes that the interpolation is performed with respect to the coordinates in the spline parametric domain $\hat{\Omega}$. We remark that in practice, the quadrature nodes in the parametric domain $\hat{\Omega}$ are employed for the interpolation. We will briefly review the EIM procedure in Section \ref{sec:eim_alg}. The transformed parameterized bilinear form $\alpha(\cdot,\cdot;\bm{\mu})$ and the linear functional $f(\cdot;\bm{\mu})$ of Equation \eqref{eq25} admit then an affinely decomposable approximation that is:
\begin{equation}\label{eq28}
\begin{aligned}
\alpha_M (u_h, v_h ; \bm{\mu})  &= \int_{\hat{\Omega}} \nabla \hat{u}_h \ g_{\alpha}^M(\boldsymbol{\xi};\bm{\mu}) \cdot {\nabla} \hat{v}_h\ \textrm{d}\hat\Omega= \sum_{{m_1}=1}^{M_{\alpha}} \theta_{m_1}^{\alpha}(\bm{\mu})\int_{\hat{\Omega}} \nabla \hat{u}_h \ \phi_{m_1}^{\alpha}(\boldsymbol{\xi})
 \cdot {\nabla} \hat{v}_h\ \textrm{d}\hat\Omega, \\
f_M (v_h ; \bm{\mu})  &= \int_{\hat{\Omega}}  \ g_f^M(\boldsymbol{\xi};\bm{\mu}) \ \hat{v}_h\ \textrm{d}\hat\Omega
= \sum_{{m_2}=1}^{M_f} \theta_{m_2}^{f}(\bm{\mu})\int_{\hat{\Omega}} \ \phi_{m_2}^{f}(\boldsymbol{\xi}) \ \hat{v}_h\ \textrm{d}\hat\Omega.
\end{aligned}
\end{equation}

We observe that the above integrals are now parameter-independent and can be computed and stored once and for all in the offline phase. The reduced basis approximation of the problem considering the EIM approximation reads: find $u_N^M \in V_N$ such that
\begin{equation}\label{eq29}
\alpha_M(u_N^M,v_N;\bm{\mu}) = f_M(v_N;\bm{\mu}), \qquad \forall v_N \in V_N.
\end{equation}

\subsection{EIM procedure}\label{sec:eim_alg}
The EIM seeks a set of basis functions and interpolation points (known as \emph{magic points}) that allow the affine expansion in Equation \eqref{eq27}. The construction of basis functions $\displaystyle \{\phi_{m_1}^{\alpha}\}_{{m_1}=1}^{M_{\alpha}}$,$\{\phi_{{m_2}}^{f}\}_{{m_2}=1}^{M_f}$ and choice of interpolation points are performed in the offline phase based on a Greedy algorithm similar to Section \ref{sec:greedy}. In the following we will briefly review the main idea behind the EIM.

Let us first consider a sufficiently large training set of parameter samples $\mathcal{P}_{train}^{EIM} \subset \mathcal{P}$. EIM selects with Greedy a set of parameter points $S_M^{\alpha}=\{\bm{\mu}_1^{\alpha},....,{\bm{\mu}}_{M_{\alpha}}^{\alpha}\}$ for the bilinear form and $S_M^{f}=\{\bm{\mu}_1^{f},....,{\bm{\mu}}_{M_f}^{f}\}$ for the linear functional. In particular, the Greedy algorithm selects at each iteration the function snapshot that is worst approximated by the current interpolant such that 
\begin{equation}\label{eq31}
\begin{aligned}
\bm{\mu}_{m+1}^{\alpha} &= \underset{\boldsymbol{\mu}\in\mathcal{P}_{train}^{EIM}}{\arg\max} \norm{g_{\alpha}(\cdot ; \bm{\mu}) - \mathcal{I}_{m}^{\boldsymbol{\xi}}g_{\alpha}(\cdot;\bm{\mu})}_{L^{\infty}(\hat{\Omega})},\\
\bm{\mu}_{m+1}^f &= \underset{\boldsymbol{\mu}\in\mathcal{P}_{train}^{EIM}}{\arg\max} \norm{g_{f}(\cdot ; \bm{\mu}) - \mathcal{I}_{m}^{\boldsymbol{\xi}}g_{f}(\cdot;\bm{\mu})}_{L^{\infty}(\hat{\Omega})}.
\end{aligned}
\end{equation}
The algorithm for the construction of the basis functions $\displaystyle \{\phi_{m_1}^{\alpha}\}_{{m_1}=1}^{M_{\alpha}}$,$\{\phi_{{m_2}}^{f}\}_{{m_2}=1}^{M_f}$ and magic points $T_M^{\alpha} = \{\boldsymbol{t}_1^{\alpha},....,{\boldsymbol{t}}_{M_{\alpha}}^{\alpha}\} \in \hat{\Omega}$, $T_M^{f} = \{\boldsymbol{t}_1^{f},....,{\boldsymbol{t}}_{M_f}^{f}\} \in \hat{\Omega}$ is presented in detail in Chapter 10.1.3 of \cite{QMN_RBspringer}.
The basis functions and magic points are selected such that for any $\bm{\mu} \in \mathcal{P}$ they fulfill the interpolation constraints 
\begin{equation}\label{eq30}
\begin{aligned}
g_{\alpha}(\boldsymbol{t}_q^{\alpha}{;\bm{\mu})} &= \mathcal{I}_M^{\boldsymbol{\xi}} g_{\alpha}(\boldsymbol{t}_q^{\alpha} ;
\bm{\mu}) 
= \sum_{{m_1}=1}^{M_{\alpha}} \theta_{m_1}^{\alpha}(\bm{\mu})\phi_{m_1}(\boldsymbol{t}_q^{\alpha}),  \ \ \ \forall  q = 1,...,M_{\alpha}, \\
g_{f}(\boldsymbol{t}_q^{f}{;\bm{\mu})} &= \mathcal{I}_M^{\boldsymbol{\xi}} g_{f}(\boldsymbol{t}_q^{f} ;
\bm{\mu}) 
= \sum_{{m_2}=1}^{M_f} \theta_{m_2}^f(\bm{\mu})\phi_{m_2}(\boldsymbol{t}_q^{f}),  \ \ \ \forall  q = 1,...,M_f.
\end{aligned}
\end{equation}
The basis functions and magic points are selected in the offline phase. In the online phase, the solution of Equation \eqref{eq30} for given values of the parameters yields the functions $\theta^{\alpha}, \theta^f$ and thus the affine expansions $g_{\alpha}^M, g_f^M$ defined in Equation \eqref{eq27}. A detailed overview of the EIM algorithm and its application to nonaffine problems is provided in \cite{Barrault2004,Maday2009,Eftang2010}.

\section{Domain decomposition with SCRBE method}\label{sec:scrbe}
We are mainly interested in parameterized geometries represented by multiple patches, which are subdomains with possibly different geometrical or physical parameters. The standard RB framework is known to have its limitations when it comes to increasing the parameter dimensions. In particular, besides the increased offline cost associated to a sufficiently large training set, the number of affine terms grows drastically with increasing dimension of the parameter space $\mathcal{P}$. This motivates domain decomposition to reduce the cost of multi-dimensional problems. The Static Condensation Reduced Basis Element (SCRBE) method is a component-based approach that facilitates efficient dimension reduction at component interiors and interfaces \cite{Eftang2013}. For our study in the context of isogeometric analysis, it seems convenient to assume components as spline patches, but of course other choices could be also possible. 

Let us now consider the union of $n_p$ patches that form the physical domain as introduced in Equation \eqref{eq9}. The boundary of each patch $\partial{\Omega}_k, k=1,...,n_p$ consists of disjoint interfaces referred to as \emph{ports} in the context of SCRBE.
In the following, we will consider only common interfaces between patches as ports, that is $\Gamma_{k,l}=\partial{\Omega_k}\cap\partial{\Omega_l} \ne \varnothing$. We define the local ports of the $k$-th patch as $\gamma_{k,j}, j=1,...,n_{\gamma}^{k}$, where $n_{\gamma}^{k}$ is the number of local ports. Now we introduce the global ports $\Gamma_p$, where $p=1,...,n_{\Gamma}$ and $n_{\Gamma}$ is the number of global ports. The connectivity of two local ports $\gamma_{k,j},\gamma_{l,j'}$ is defined through the index set $\pi_p = \{(k,j),(l,j')\}$. Moreover, we also introduce a local to global port index mapping $\mathcal{G}_k$ such that $\mathcal{G}_k(j)=p$ for $p$ such that $(k,j) \in \pi_p$. 

The static condensation eliminates the degrees of freedom in the interior of the patch and expresses them in terms of degrees of freedom at the interface, which can be further reduced to only a few port modes. For the purpose of formulating the static condensation procedure, let us first introduce the \emph{bubble} space of the $k$-th patch as
\begin{equation}\label{eq32}
V_{h,b}^{k} = \{ v_h  \in V_h^{k} : v_h \vert_{\gamma_{k,j}} = 0, \ j=1,...,n_{\gamma}^{k}\},
\end{equation}
where $V_h^{k}$ is the finite dimensional subspace associated to the $k$-th patch. We can now define the bubble space in the physical domain as $V_{h,b}=\oplus_{k=1}^{n_p} V_{h,b}^{k}$. Furthermore, we define the port space by restricting $V_{h,b}^{k}$ to $\gamma_{k,j}$, $j=1,...,n_{\gamma}^{k}$ as 
\begin{equation}\label{eq32b}
V_{h,\gamma}^{k,j} = \{ v_h  \in V_h^{k} \vert \ v_h \vert_{\gamma_{k,j}} \ne 0\} = \text{span}\{\chi_1^{k,j},...,\chi_{\mathcal{N}^{\gamma}_{k,j}}^{k,j}\}.
\end{equation}
The dimension of the local port space is $ \mathcal{N}^{\gamma}_{k,j}=\dim(V_{h,\gamma}^{k,j})$. Since conforming port spaces are considered, we also define the dimension of the global port space associated to the $p$-th global port as $\mathcal{N}^{\Gamma}_{p}=\mathcal{N}^{\gamma}_{k,j}$. Now let us consider the subdomain $\Omega_k(\bm{\mu}^{(k)})$ associated to the $k$-th patch and introduce the functions
\begin{equation}\label{eq33}
\phi_h^{k,j,r}(\bm{\mu}^{(k)}) = b_h^{k,j,r}(\bm{\mu}^{(k)}) + \psi_r^{k,j},
\end{equation}
where the functions $\psi_r^{k,j} \in V_{h}^{k} $ are harmonic extensions of the basis functions $\chi_r^{k,j} \in V_{h,\gamma}^{k,j}$ for $r=1,.., \mathcal{N}^{\gamma}_{k,j}$. The bubble functions $b_h^{k,j,r}(\bm{\mu}^{(k)}) \in V_{h,b}^{k}$ satisfy
\begin{equation}\label{eq34}
\alpha_k(b_h^{k,j,r}(\bm{\mu}^{(k)}),v;\bm{\mu}^{(k)}) = - \alpha_k(\psi_r^{k,j},v;\bm{\mu}^{(k)}), \qquad \forall v \in V_{h,b}^{k},
\end{equation}
where $j=1,..,n_{\gamma}^{k}$ and $r=1,.., \mathcal{N}^{\gamma}_{k,j}$. Moreover, we define the bubble functions $b_{h,f}^{k}(\bm{\mu}^{(k)}) \in V_{h,b}^{k}$ associated to the right-hand side of the $k$-th patch such that
\begin{equation}\label{eq35}
\alpha_k(b_{h,f}^{k}(\bm{\mu}^{(k)}),v;\bm{\mu}^{(k)}) = f_k(v;\bm{\mu}^{(k)}), \qquad \forall v \in V_{h,b}^{k}.
\end{equation}
Now we can express the solution on each patch as
\begin{equation}\label{eq36}
u_h(\bm{\mu}) \vert_{\Omega_k} = b_{h,f}^{k}(\bm{\mu}^{(k)}) + \sum_{j=1}^{n_{\gamma}^{k}} \sum_{r=1}^{\mathcal{N}^{\gamma}_{k,j}} \hat{u}_{\mathcal{G}_k(j),r}(\bm{\mu})\phi_h^{k,j,r}(\bm{\mu}^{(k)}),
\end{equation}
where the solution coefficients $\hat{u}_{\mathcal{G}_k(j),r}(\bm{\mu})$ are unknowns to be determined in the following, with $r=1,..,\mathcal{N}^{\Gamma}_{p}$ and $p=1,...,n_{\Gamma}$. We remark that the functions $\phi_h^{k,j,r}$ are glued together at global ports, that is for global port $\pi_p = \{(k,j),(l,j')\}$ we define $\Phi_{p,r}=\phi_h^{k,j,r} + \phi_h^{l,j',r}$ and assume that functions are extended by zero outside their domain of definition. The skeleton space can be then defined as
\begin{equation}\label{eq37}
V_h^{\mathcal{S}} = \text{span}\{\Phi_{p,r}(\bm{\mu}), \ r=1,.., \mathcal{N}^{\Gamma}_{p}, \ p=1,...,n_{\Gamma}\} \subset V_h.
\end{equation}
Let us now express the global solution as
\begin{equation}\label{eq38}
u_h(\bm{\mu}) = \sum_{k=1}^{n_p} b_{h,f}^{k}(\bm{\mu}^{(k)}) + \sum_{p=1}^{n_{\Gamma}} \sum_{r=1}^{\mathcal{N}^{\Gamma}_{p}} \hat{u}_{p,r}(\bm{\mu}) \Phi_{p,r}(\bm{\mu}).
\end{equation}
To complete the static condensation, we insert the solution representation \eqref{eq38} into Eq. \eqref{eq2} and restrict the test space to $V_h^{\mathcal{S}}$, that is
\begin{equation}\label{eq39}
\sum_{p=1}^{n_{\Gamma}} \sum_{r=1}^{\mathcal{N}^{\Gamma}_{p}} \hat{u}_{p,r}(\bm{\mu}) \alpha(\Phi_{p,r}(\bm{\mu}), v ;\bm{\mu}) = f(v;,\bm{\mu}) - \sum_{k=1}^{n_p} \alpha(b_{h,f}^k(\bm{\mu}^{(k)}),v;\bm{\mu}), \quad \forall v \in V_h^{\mathcal{S}}.
\end{equation}
There are two stages of reduction that we will consider in the following: a reduced basis approximation in the interior of each patch, referred to as bubble approximation, and a port reduction that retains only a few dominant port modes at interfaces. For this purpose, we introduce the reduced bubble space $V_{N_b} = \oplus_{k=1}^{n_p} V_{N_b}^{(k)} \subset V_{h,b}$
where $V_{N_b}^{(k)}$ is the reduced bubble space associated to the $k$-th patch. Moreover, we introduce a reduced port space $V_{N,\gamma}^{k,j}=\text{span}\{\chi_1^{k,j},...,\chi_{{n}^{\gamma}_{k,j}}^{k,j}\} \subseteq V_{h,\gamma}^{k,j}$ associated to the $k$-th patch and $j$-th interface, where ${n}^{\gamma}_{k,j} \leq \mathcal{N}^{\gamma}_{k,j}$ is the dimension of the reduced port space. The reduced bubble and port spaces can be constructed in a similar fashion to the techniques described in Section \ref{sec:basis}. More details on their construction and properties can be found in \cite{Eftang2013,Smetana2016}. In this work, we employ the Greedy algorithm to construct the reduced bubble space and we use the POD to obtain the reduced port space. Thus, we can obtain reduced basis approximations of the bubble functions $\tilde{b}^{k,j,r}(\bm{\mu}^{(k)}) \in V_{N_b}^{(k)}$ for $r=1,...,{n}^{\gamma}_{k,j}$ and $\tilde{b}_{f}^{k}(\bm{\mu}^{(k)}) \in V_{N_b}^{(k)}$ in Eqs. \eqref{eq34}-\eqref{eq35}.

 It is worthwile noting that if ports are not mutually disjoint, cross-points or wire-baskets can be treated as discussed in \cite{Huynh2013a,Hetmaniuk2010}. The bilinear forms and linear functionals in Equations \eqref{eq34},\eqref{eq35} might entail a nonaffine parametric dependence. To recover an affine decomposition in the form of Equation \eqref{eq14} one can employ the EIM as discussed in Section \ref{sec:eim}. In the following, we will illustrate the presented procedure on a numerical example of a parameterized multipatch geometry.

\section{Numerical Results}\label{sec:example}
In this section we present some numerical results to assess the performance of the strategy developed in this work. For the purpose of illustrating the presented theoretical framework, we consider a curved multipatch geometry with a multi-dimensional geometrical parameterization. The aim of this numerical study is to investigate the capabilities of splines in combination with the presented reduced basis framework and demonstrate the feasibility of this approach on a simple test case. The full order solution snapshots have been computed using spline approximations with the open-source Octave/Matlab isogeometric package \emph{GeoPDEs} \cite{Vazquez2016}. The reduced basis approximations have been obtained with the open-source package \emph{redbKIT} \cite{QMN_RBspringer,redbKIT} as an external library. 

\begin{figure}[!htb]
	\centering
	\includegraphics[width=1.0\textwidth]{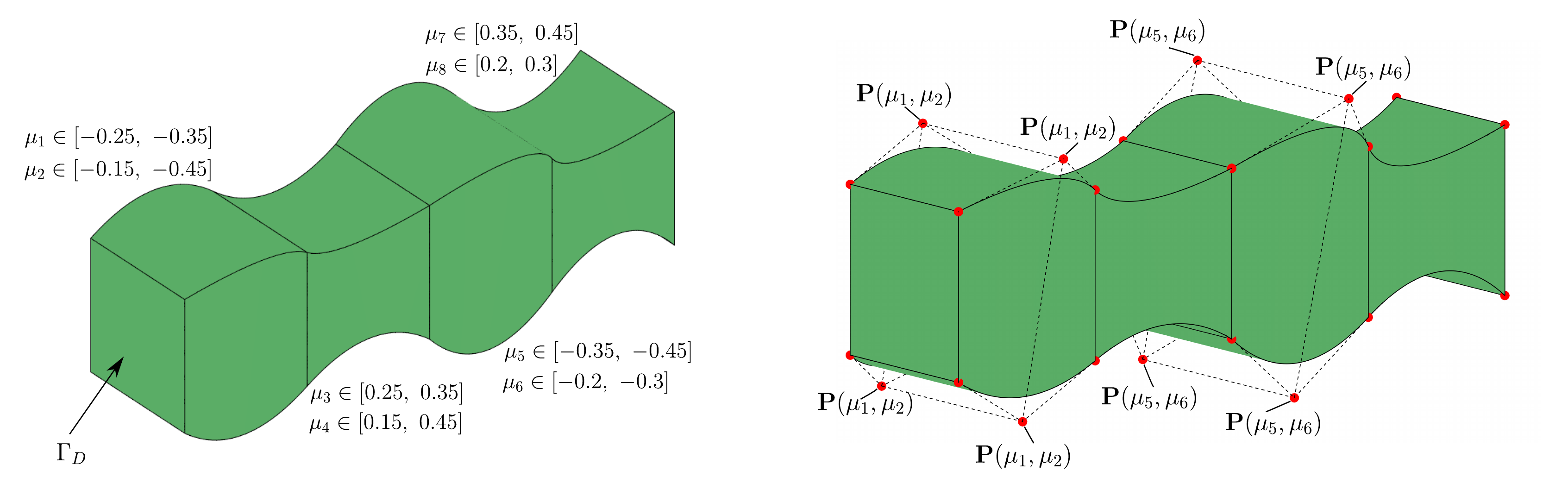} \\
	\caption{Geometry and parameterization.}\label{fig:system}
\end{figure}

Let us consider the geometry depicted in Figure \ref{fig:system}. The domain is parameterized using eight geometrical parameters $\bm{\mu}=(\mu_1,\mu_2,\mu_3,...,\mu_8)$. The parameters associated to each patch are given as
\begin{equation}\label{eq40}
\bm{\mu}^{(1)}=[\mu_1,\mu_2], \quad  \bm{\mu}^{(2)}=[\mu_3,\mu_4], \quad \bm{\mu}^{(3)}=[\mu_5,\mu_6], \quad \bm{\mu}^{(4)}=[\mu_7,\mu_8].
\end{equation}
 We remark that configurations of the domain for certain values of the geometrical parameters may be difficult to obtain with one single patch without introducing high mesh distortions. Therefore, the domain is constructed using 4 patches, while the patch interfaces are visible in Figure \ref{fig:system}. The geometrical parameters are defined on a coarse geometry by parameterizing the control points as described in Section \ref{sec:iga}. Then the geometry is refined to obtain the solution snapshots. The parameterized control points are illustrated in Figure \ref{fig:system} for patches 1 and 3, while the same control points are also parameterized for patches 2 and 4. The polynomial degree is set to $p=3$ to represent the geometry and approximate the solution. We construct a full order model using a mesh with 1225 degrees of freedom and 64 elements per patch. Homogeneous Dirichlet boundary conditions are applied on the boundary $\Gamma_D$ and the source term is chosen as $\tilde{f}=2xyz$ for all patches.

\begin{figure}[!htb]
	\centering
	\includegraphics[height=0.5\textwidth]{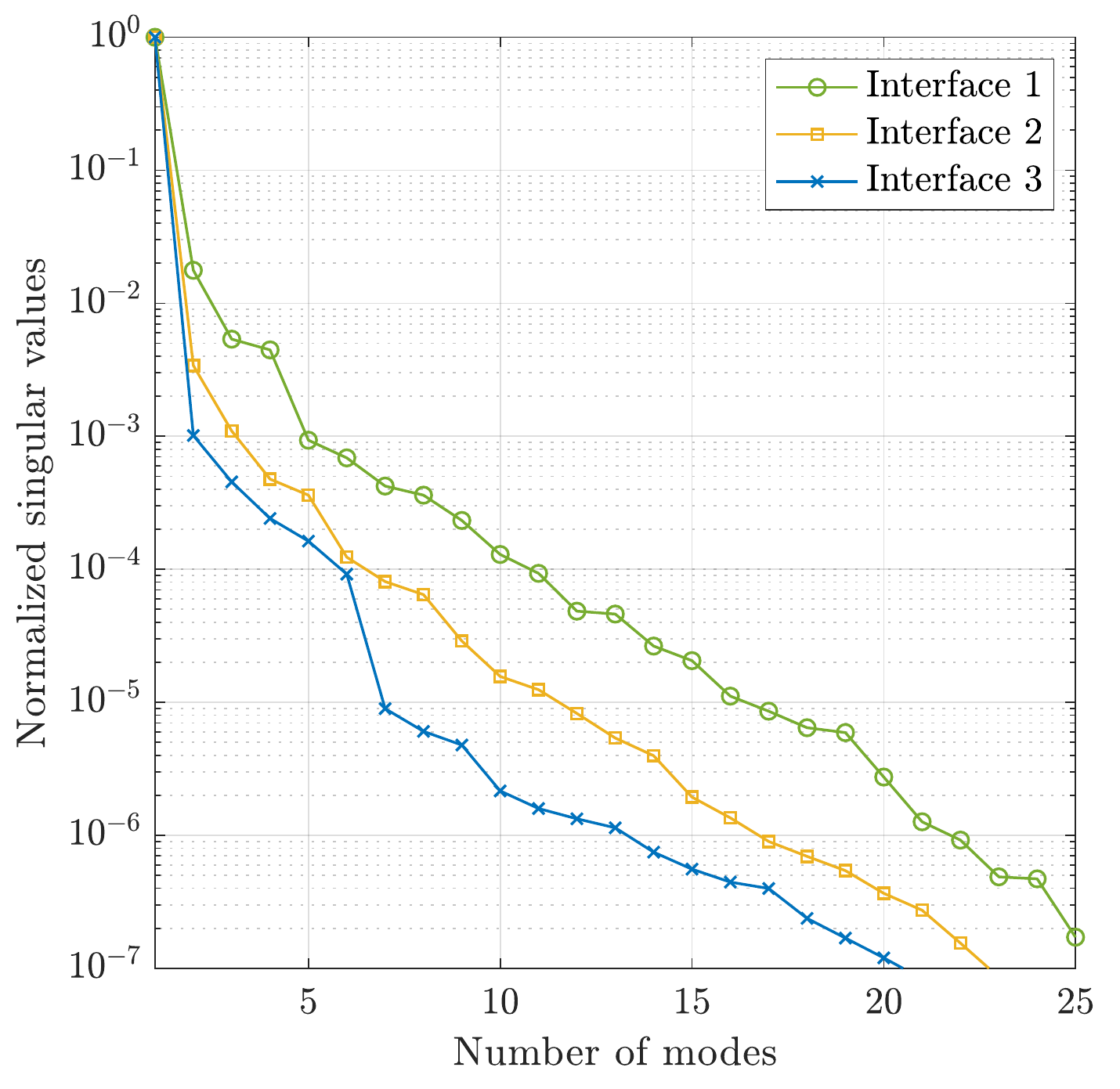} \\
	\caption{Singular value decay vs. number of port modes for each interface.}\label{fig:port}
\end{figure}

\newpage Let us now construct a reduced basis approximation based on the concepts introduced in Sections \ref{sec:rb}-\ref{sec:scrbe}. Since the problem at hand exhibits a multi-dimensional parameterization, we exploit the SCRBE procedure to reduce the dimensionality of the parameter space. The first step towards an efficient offline/online decomposition is to construct a reduced port space spanned by a set of port modes. For the sake of simplicity, we construct empirical port modes by applying the POD on a set of snapshot solutions computed at the interfaces between patches. These are pre-computed in the offline phase by employing a NUBRS multipatch approximation as full order model. We remark that this serves the purpose of our numerical study, while in principle other port mode constructions are possible \cite{Eftang2013b,Smetana2016}. To construct the snapshot solutions we consider a set of random parameter values of dimension $N_s = 25$, which are obtained by Latin hypercube sampling \cite{McKay1979,Cochran2007}. The singular value decay of the POD depicted is in Figure \ref{fig:port} for each interface. We observe a rapid decay, which indicates that the dimension of the problem at the interface can be effectively reduced to a small number of port modes ($\leq 25$).

\begin{figure}[!h]
     \begin{subfigure}[b]{0.49\textwidth}
         \centering
         \includegraphics[width=\textwidth]{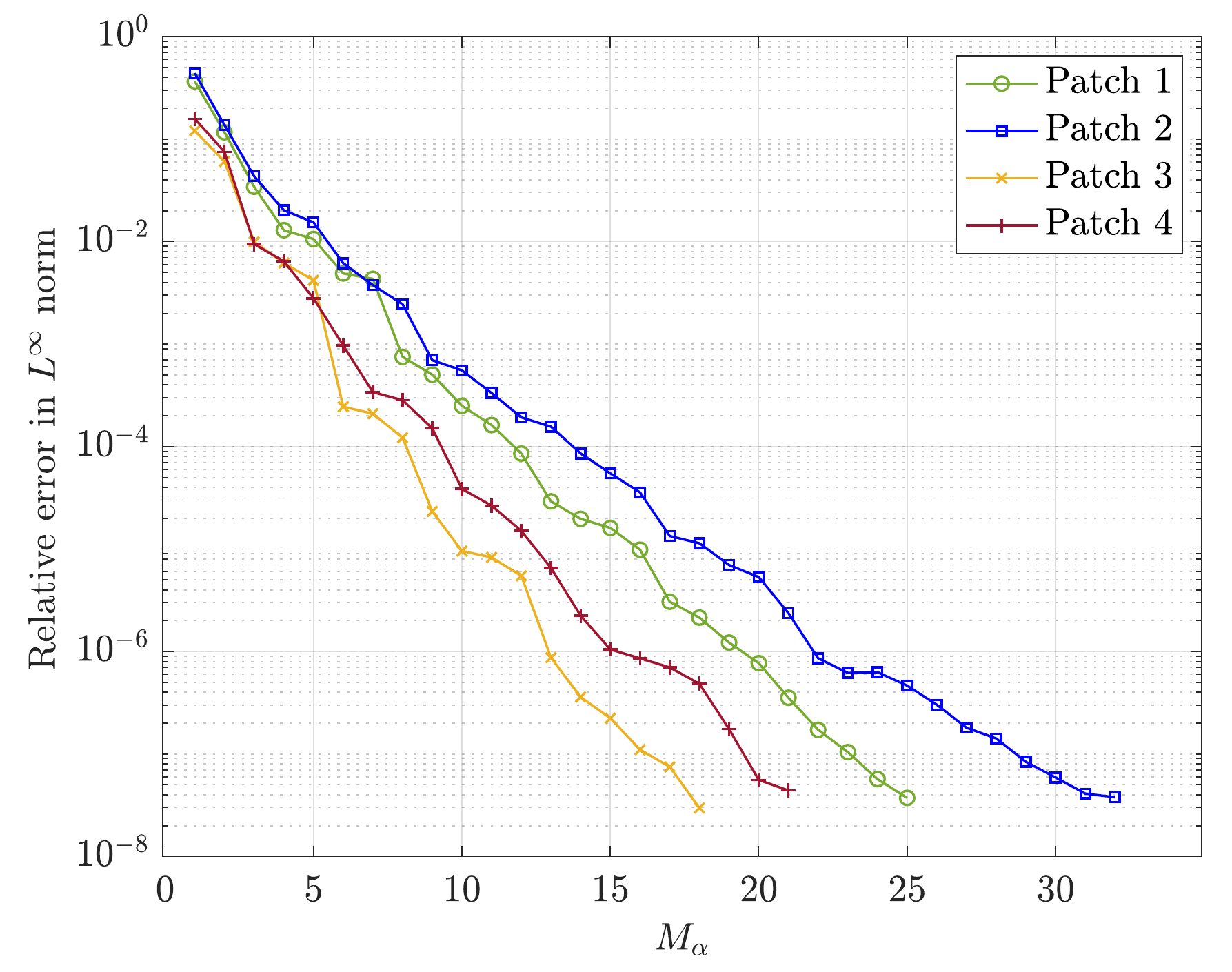}
         \caption{Stiffness matrix}
         \label{fig:EIM_mat}
     \end{subfigure}
     \begin{subfigure}[b]{0.49\textwidth}
         \centering
         \includegraphics[width=\textwidth]{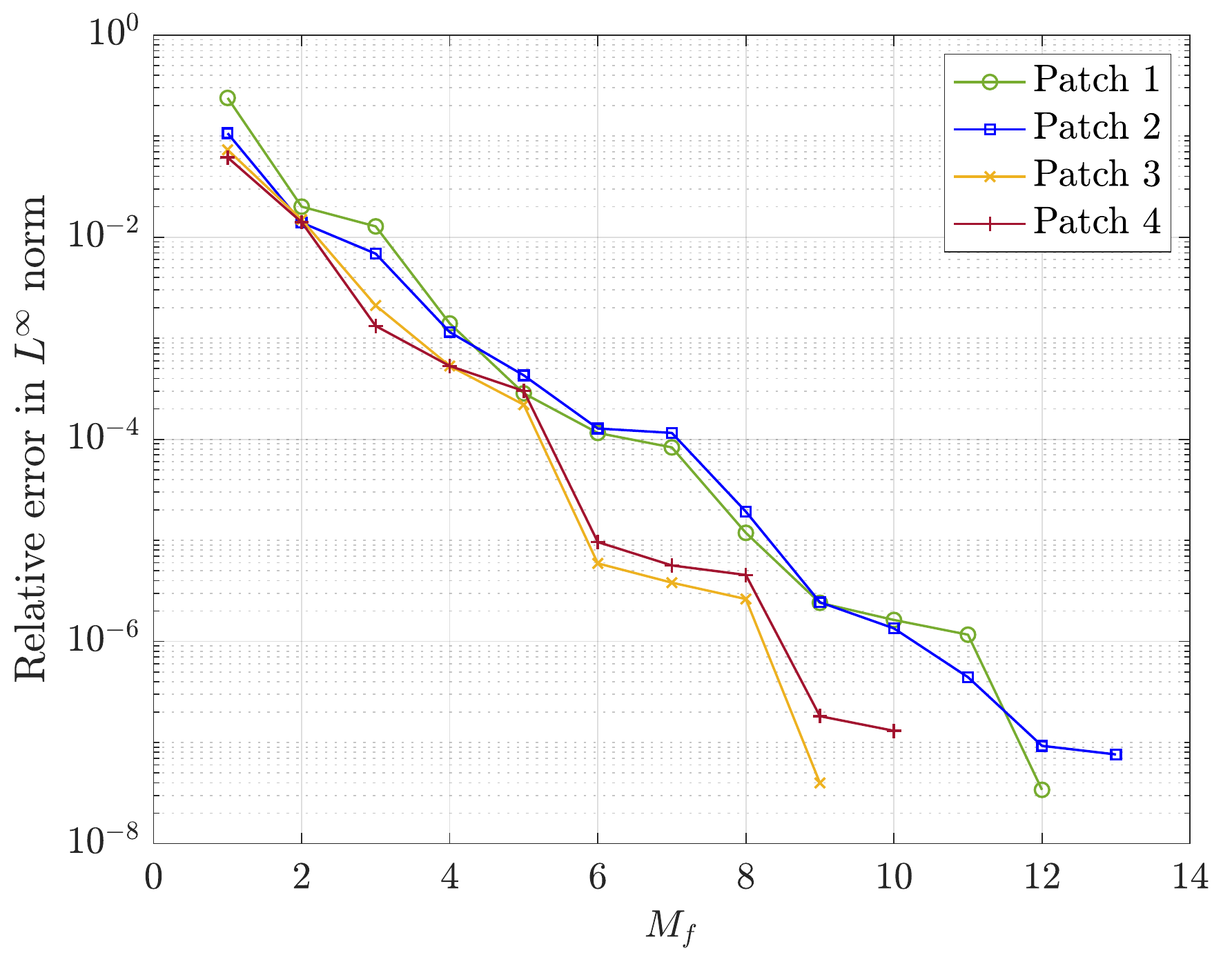}
         \caption{Right-hand side}
         \label{fig:EIM_rhs}
     \end{subfigure}
        \caption{Error decay in $L^\infty$ norm of the EIM approximations.}
      \label{fig:EIM}
\end{figure}

The next step in the offline phase is to construct the EIM approximations. The prescribed tolerance of the Greedy algorithm is set to $\varepsilon = 10^{-7}$ and the dimension of the parameter set employed for the training is $\abs{\mathcal{P}_{train}^{EIM}}=250$. The obtained number of basis functions $M_{\alpha}$ and $M_f$ for the affine approximation of the matrix and right-hand side, respectively, are summarized in Table \ref{tab:cost}. Moreover, Figure \ref{fig:EIM} depicts the error decay over the number of basis functions for each patch. The error analysis is performed on a test sample set of dimension $N_{test}^{EIM}=100$ that consists of random parameter values obtained by Latin hypercube sampling. We observe that a small number of basis functions ($M_{\alpha} \le 33$, $M_f \le 14$) is sufficient to achieve an accuracy of $10^{-7}$.

\begin{figure}[!htb]
     \begin{subfigure}[b]{0.49\textwidth}
         \centering
         \includegraphics[width=\textwidth]{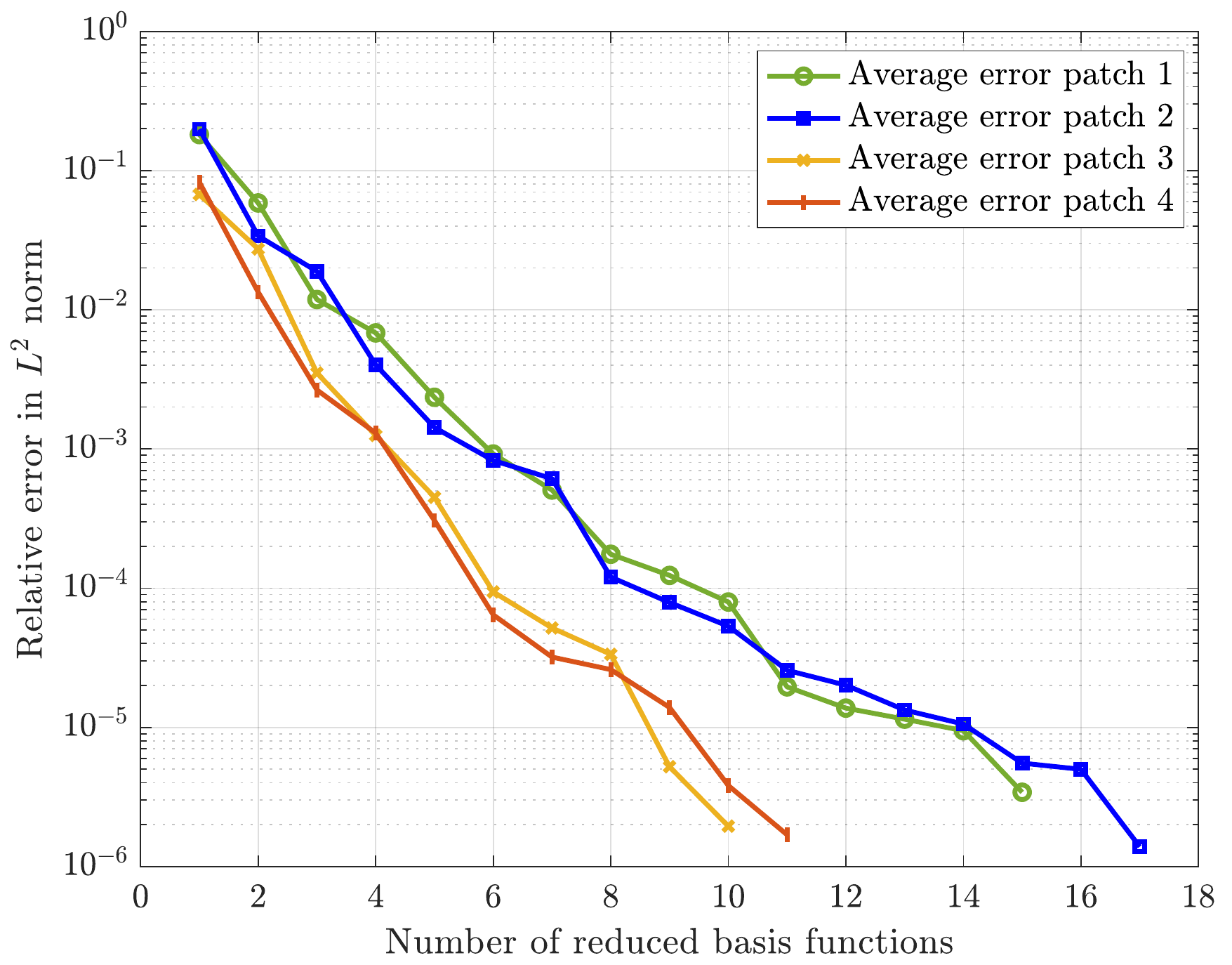}
         \caption{Average error}
         \label{fig:ROM_av}
     \end{subfigure}
     \begin{subfigure}[b]{0.49\textwidth}
         \centering
         \includegraphics[width=\textwidth]{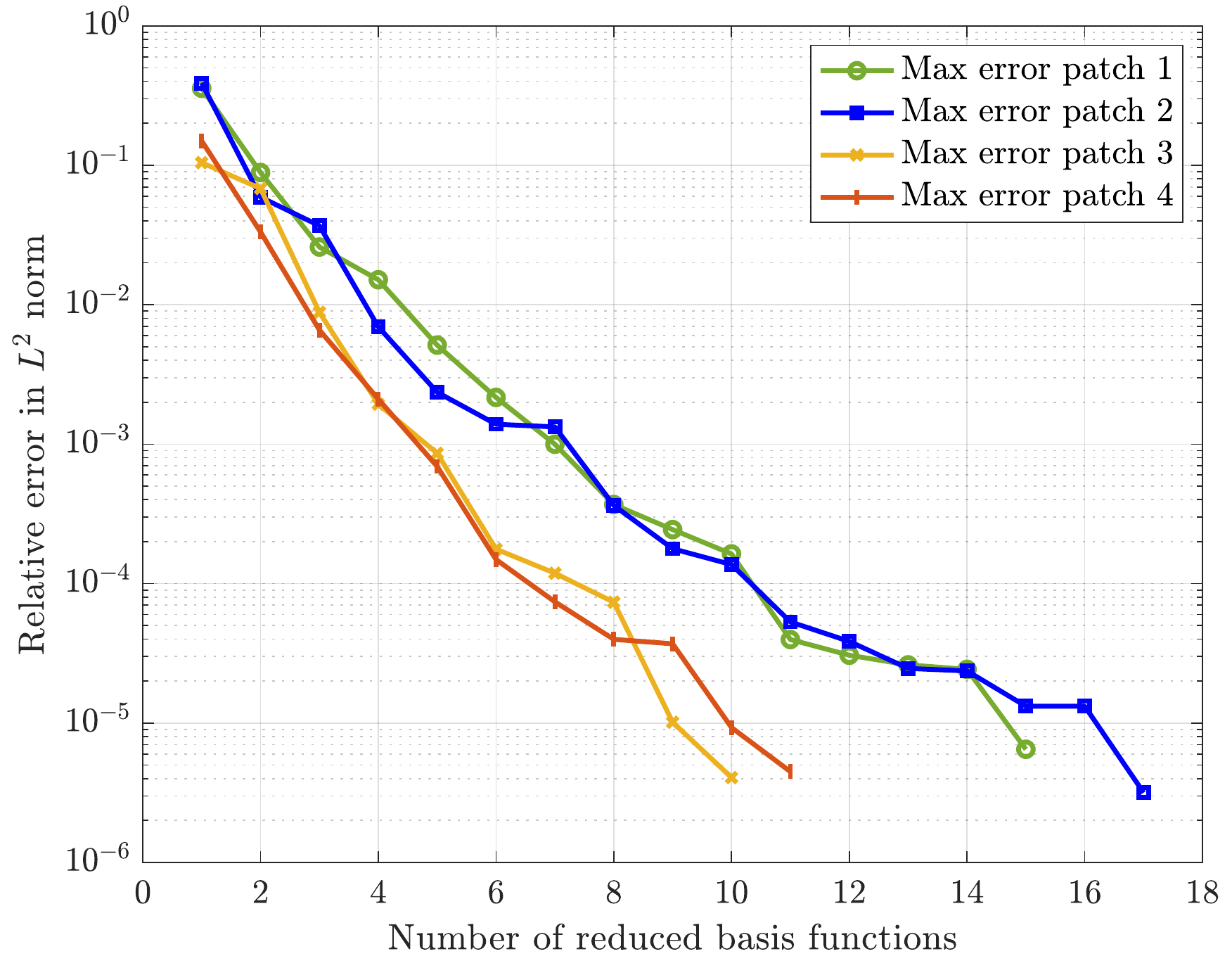}
         \caption{Max. error}
         \label{fig:ROM_max}
     \end{subfigure}
        \caption{Error decay of the reduced basis bubble approximation for each patch.}
      \label{fig:ROM}
\end{figure}

\begin{figure}[!htb]
     \begin{subfigure}[b]{0.49\textwidth}
         \centering
         \includegraphics[width=\textwidth]{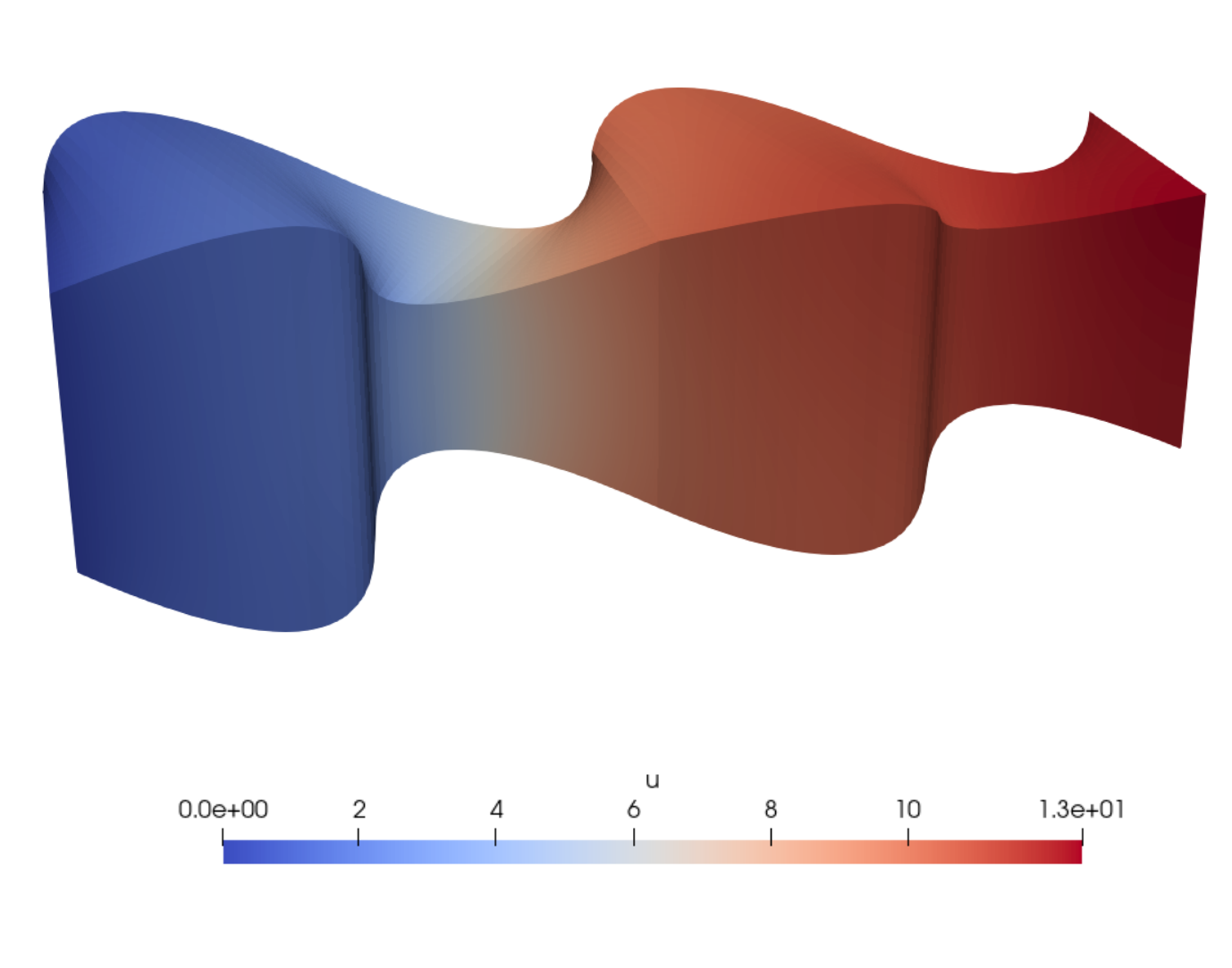}
         \caption{Full order model}
         \label{fig:plot_FOM}
     \end{subfigure}
     \begin{subfigure}[b]{0.49\textwidth}
         \centering
         \includegraphics[width=\textwidth]{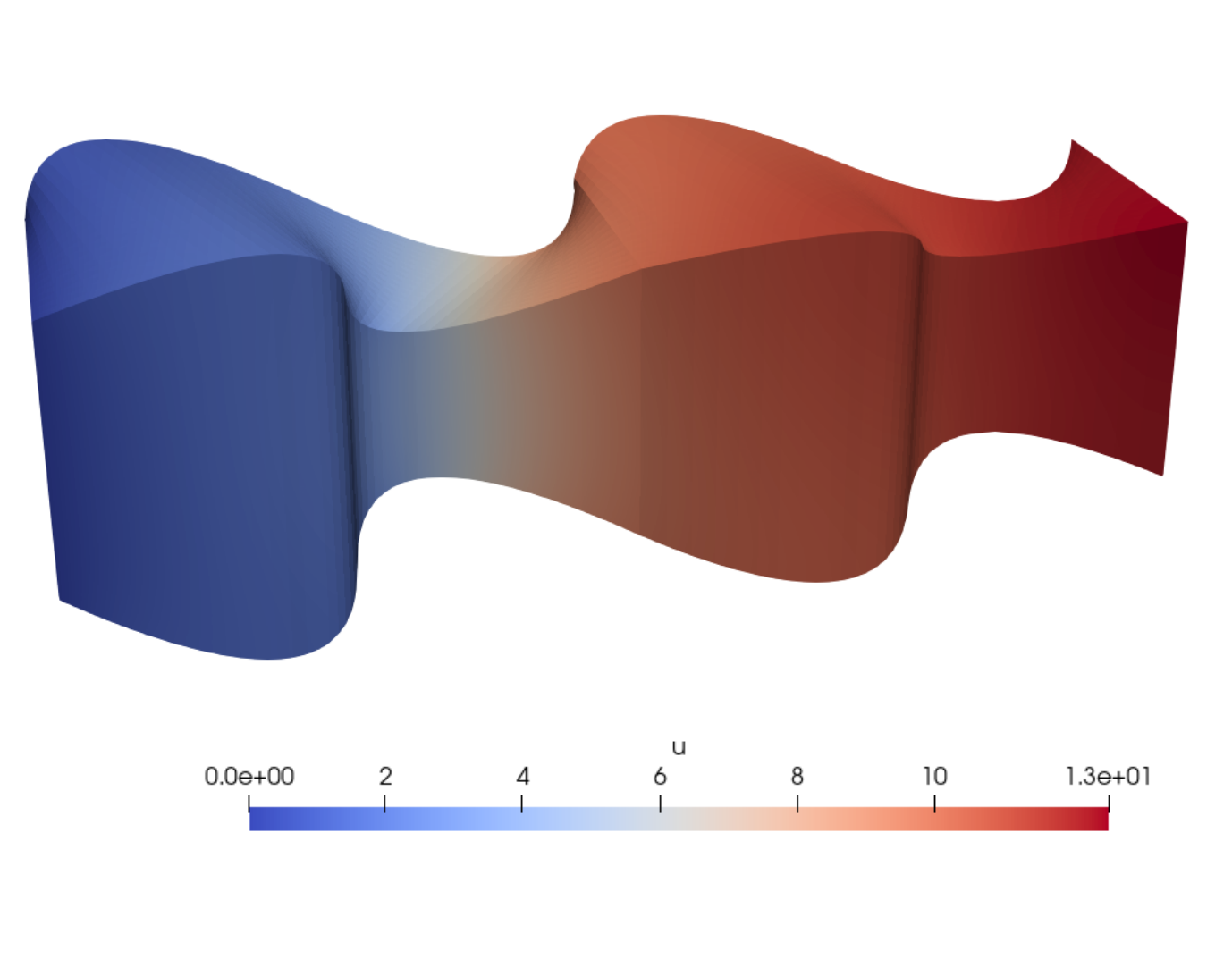}
         \caption{Reduced order model}
         \label{fig:plot_ROM}
     \end{subfigure}
        \caption{Comparison of the solution between the full and reduced order model.}
      \label{fig:plot}
\end{figure}

We now construct reduced basis bubble approximations for each patch with SCRBE as described in Section \ref{sec:scrbe}. The prescribed tolerance for the Greedy algorithm is set to $\varepsilon = 10^{-5}$ and the dimension of the parameter set employed for the training is $\abs{\mathcal{P}_{train}}=250$. The error decay in the $L^2$-norm for the reduced bubble functions associated to the right-hand side (cf. \eqref{eq35}) is depicted in Figure \ref{fig:ROM} . The error analysis is performed on a test sample set of dimension $N_{test}=30$ and random parameter values. Similarly to the EIM approximations, a small set of basis functions ($\leq 17)$ is sufficient to achieve an accuracy of $10^{-5}$. Since the EIM requires more basis functions to approximate the stiffness matrix of patches 1 and 2, the number of bubble basis functions is also higher compared to patches 3 and 4. Table \ref{tab:cost} summarizes the obtained number of bubble basis functions and the computational cost for each patch. An average online evaluation of the bubble RB approximation per patch takes $42$ ms. Regarding the reduced bubble functions corresponding to the interface (cf. \eqref{eq34}), the number of reduced basis functions associated to each port mode is on average 16. The training can be performed in parallel to speed up the computation in the offline phase. On average, the offline training time for Greedy is $1.61$ min and the online evaluation is $29.3$ ms for each port mode. Figure \ref{fig:plot} shows a qualitative comparison of the solution between the full and reduced order model.

 \begin{table}[!htb]
	\centering
	\caption{Number of bubble basis functions and computational cost for each patch.}\label{tab:cost}
	\begin{tabular}{cccccc} \hline 
		 & $M_{\alpha}$ &  $M_f$ & $N_{b}$ & Offline Greedy time [min] & Online CPU time [ms] \\\hline
		Patch 1  & 28 & 13 & 17 & 2.16 & 44.2 \\ 
	    Patch 2  & 33 & 14 & 16 & 1.72 & 41.2 \\ 
		Patch 3  & 19 & 10 & 10 & 1.01 & 46.4 \\ 
	    Patch 4 & 23 & 11 & 11 & 1.70 & 36.8 \\\hline
	\end{tabular}
 \end{table}

\newpage

\section{Conclusions}
In this work we presented a reduced basis framework for nonaffine problems approximated by splines. We employed EIM and formulated the reduced basis approximation by exploiting the exact geometric map. To tackle problems that exhibit a high dimensional parameter space, we considered domain decomposition with the SCRBE method. The procedure was illustrated on a 3D example with a multi-dimensional geometrical parameterization. The obtained results aim to investigate the efficiency of the developed procedure. Indeed, the dimension of the obtained reduced basis is small and the online evaluation is rapid. The application of the presented framework to more involved problems with complex geometries of industrial relevance is a subject of future work. Finally, another interesting direction for future research is the integration into a parametric optimization process in order to exploit the potential of the reduced basis framework.

\section*{Acknowledgments}
The financial support of the Swiss Innovation Agency (Innosuisse) under Grant No. 46684.1 IP-EE is gratefully acknowledged. We also wish to thank Dr. David Knezevic and Dr. Jonas Ballani from Akselos S.A. for fruitful discussions.

\bibliography{mybibfile}

\end{document}